\documentclass[11pt]{article} \makeatother
\usepackage{amsfonts,amsmath,amssymb}
\allowdisplaybreaks[2]

\begin{document}

\numberwithin{equation}{section}
\newtheorem{theorem}{\ \ \ \ Theorem}[section]
\newtheorem{proposition}[theorem]{\ \ \ \ Proposition}
\newtheorem{lemma}[theorem]{\ \ \ \ Lemma}
\newtheorem{remark}{\ \ \ \ Remark}
\newcommand{\be}{\begin{equation}}
\newcommand{\ee}{\end{equation}}
\newcommand\bes{\begin{eqnarray}} \newcommand\ees{\end{eqnarray}}
\newcommand{\bess}{\begin{eqnarray*}}
\newcommand{\eess}{\end{eqnarray*}}
\newcommand\ds{\displaystyle}

 \begin{center} {\bf\Large General decay of the solution for a viscoelastic wave  equation }\\[2mm]
  {\bf\Large  with a time-varying delay term in the internal feedback}
\\[4mm]
  {\large Wenjun Liu}\\[1mm]
 College of Mathematics and Statistics, Nanjing University of
Information Science and Technology, Nanjing 210044, China \\[1mm] E-mail: wjliu@nuist.edu.cn \\[2mm]
\end{center}

\setlength{\baselineskip}{17pt}{\setlength\arraycolsep{2pt}

\begin{quote}
\noindent {\bf Abstract:} In this paper we consider  a viscoelastic wave  equation
  with a time-varying delay term, the coefficient of which is not necessarily positive. By introducing suitable energy and Lyapunov functionals, under suitable assumptions, we establish a general energy decay result from which
the exponential and polynomial types of decay are only special cases.


\end{quote}

\setlength{\baselineskip}{17pt}{\setlength\arraycolsep{2pt}

\section{Introduction}

In this work, we investigate the   following viscoelastic wave
equation with a linear damping and a time-varying delay term in the internal feedback
 \bes\left\{\begin{array}{ll}
\displaystyle  u_{tt}(x,t)-\Delta u(x,t)+\int_{0}^{t}g(t-s)\Delta
u(x,s)\, {\rm d}s\quad &  \medskip\\\medskip
\quad \quad \quad \quad \quad \quad +a_{0}u_{t}(x,t)+a_{1}u_{t}(x,t-\tau(t))=0,
&(x,t)\in \Omega \times (0,\infty ),\\\medskip
  \displaystyle u(x,t)=0,  &(x,t)\in \partial\Omega \times [0,\infty ),\\\medskip
\displaystyle u(x,0)=u_{0}(x),\quad u_{t}(x,0)=u_{1}(x), \quad   &x
\in  \Omega,\\
  u_{t}(x,t)=f_{0}(x,t)\quad  &(x,t)\in \Omega\times [-\tau(0),0),
\end{array}\right.\label{1.1}
 \ees
where $\Omega$
 is a bounded domain of $\mathbb{R}^n(n\geq 2)$ with a boundary $\partial\Omega$
 of class $C^2$, $a_0$ and $a_1$ are real numbers with $a_0 > 0$, $\tau(t)>0$ represents the time-varying delay, and the initial datum
$u_0, u_1, f_0$ are given functions belonging to suitable spaces.

The viscoelastic wave equation  without delay (i.e., $a_1 = 0$), has been considered by many authors during the past decades.
 Cavalcanti et al. \cite{cvs2002} studied
$$u_{tt}-\Delta u+\int_{0}^{t}g(t-\tau)\Delta u(\tau){\rm d}\tau+a(x)u_{t}+|u|^\gamma u=0, \quad (x, t)\in \Omega\times(0, \infty), $$
for $a:\Omega\rightarrow \mathbb{R^+}$, a function, which may be null on a part of the domain $\Omega$. Under the conditions that $a(x)\geq a_{0}>0$
on $\omega\subset\Omega$, with $\omega$ satisfying some geometry restrictions and $$-\xi_{1} g(t)\leq g'(t)\leq -\xi_{2} g(t), \quad t\geq0, $$
the authors established an exponential rate of decay. Berrimi and Messaoudi \cite{bm2004} improved Cavalcanti's result by introducing a different functional which allowed to weak the conditions on both $a$ and $g$.   In \cite{co2003}, Cavalcanti et al. considered
$$u_{tt}-k_{0}\Delta u+\int_{0}^{t}\mathrm{div}[a(x)g(t-\tau)\nabla u(\tau)]{\rm d}\tau+b(x)h(u_{t})+f(u)=0, $$
under similar conditions on the relaxation function $g$ and $a(x)+b(x)\geq\rho>0$, for all $x\in\Omega$. They improved the result of \cite{cvs2002} by establishing exponential stability for $g$ decaying exponentially and $h$ linear and polynomial stability for $g$ decaying polynomially and $h$ nonlinear.  Berrimi and  Messaoudi \cite{bm06} considered
$$u_{tt}-\Delta u+\int_{0}^{t}g(t-\tau)\Delta u(\tau){\rm d}\tau=|u|^{p-2}u, \quad p>2$$
in a bounded domain. They   showed, under weaker conditions than those in \cite{co2003}, that the solution is global and decay in a polynomial or exponential fashion when the initial data is small enough. Then   Messaoudi \cite{m2008} improved this result by establishing a general decay of energy which is similar to the relaxation function.
  For other related works, we refer the readers to \cite{ccfs2001,liu1,l2009,ly2011,pk2010,am9,sw2006,y2009,w2012,z1990} and the references therein.

  In recent years, the control of PDEs with time delay effects has become an active area of research, see for instance \cite{adbb1993,l2002,sxl2012,sb1980} and the references therein. The presence of delay may be a source of instability. For example, it was proved in
\cite{d1988,dlp1986,np2006,np2008,xyl2006} that an arbitrarily small delay may destabilize
a system which is uniformly asymptotically stable in the absence of delay unless
additional conditions or control terms have been used. In
\cite{np2006}, Nicaise and Pignotti examined \eqref{1.1}  with $g\equiv 0, a_0 > 0, a_1 > 0$ and $\tau(t)\equiv \tau$ be a constant delay  in the case of mixed homogeneous Dirichlet-Neumann boundary conditions, under a
geometric condition on the Neumann part of the boundary. Assuming that
$0\leq a_1<a_0,$
 a stabilization result is given, by using a suitable observability estimate
and inequalities
obtained from Carleman estimates for the wave equation due to Lasiecka et al. in \cite{lty1999}.
However, for the opposite case $a_1\ge a_0$, they were able
to construct a sequence of delays for which the corresponding solution is unstable.
The same results were obtained for the case when both the damping and the delay
act on the boundary, see also \cite{aabm2009} for the treatment of this problem in more
general abstract form. Kirane and Said-Houari \cite{ks2011} considered \eqref{1.1}  with $a_0 > 0, a_1 > 0$ and $\tau(t)\equiv \tau$ be a constant delay  in the case of initial and Dirichlet boundary conditions. They
established general energy decay results under the condition that $0\leq a_1\le a_0$.

Recently, the stability of PDEs with time-varying delays was studied in \cite{crs2007,fo2007,np2011,nvf2009,of2007}. In  \cite{nvf2009}, Nicaise et al. analyzed the exponential stability of the heat and wave equations
with time-varying boundary delay in one space dimension, under the condition $0\leq a_1<\sqrt{1-d}\, a_{0}$, where $d$ is a constant such that $\tau'(t)\leq d<1, \forall\ t>0.$  In
\cite{np2011}, Nicaise and Pignotti
studied the stabilization problem by interior damping of the wave equation with
  internal time-varying delay feedback and obtained exponential stability estimates by
introducing suitable Lyapunov functionals, under the condition $|a_1|<\sqrt{1-d}\, a_{0}$ in which the positivity of the coefficient $a_1$ is not necessary.

Motivatied by these results, we investigate in this paper system \eqref{1.1} under suitable assumptions  and prove a general decay result from which the exponential
and polynomial types of decay are only special cases.  Our main novel contribution is an extension
of previous results from \cite{ks2011,np2006} to time-varying delays with not necessarily positive coefficient $a_1$  of the delay term.
This extension is not straightforward due to the loss of translation-invariance.  For our purpose, we introduce new energy and Lyapunov functionals, which take into
account the dependence of the delay with respect to time.

The paper is organized as follows. In Section \ref{s2} we present some assumptions  and state the main result.
The general  decay result  is proved in Sections \ref{s3}.

\section{Preliminaries and main result} \label{s2}
In this section, we present some assumptions  and state the main result. We
use the standard Lebesgue space $L^2(\Omega)$ and the Sobolev space $H^1_0(\Omega)$ with their
usual scalar products and norms. Throughout this paper, $C_i$ is used to denote a
generic positive constant.

For the relaxation function $g$, we assume the following (see \cite{ly2011,m2008}):
 \begin{quote}
(G1) $g(t): [0, \infty)\to (0, \infty)$ is a non-increasing $C^1$ function such
that
 \[1-\int_{0}^{\infty }g(s)\, {\rm d}s = l>0.\]

(G2) There exists a positive non-increasing differentiable function $\xi(t)$ such
that
  $$g'(t)\leq - \xi(t)g(t),\quad t\geq 0, $$
and
 $$\int_0^{+\infty}\xi(t){\rm d}t=\infty. $$
\end{quote}

For the time-varying delay, we assume as in \cite{np2011} that there exist positive constants $\tau_{0}$,
$\overline{\tau}$ such that
\begin{align}
0<\tau_{0}\leq\tau(t)\leq\overline{\tau},\quad \forall\ t>0.
\label{1.2}
\end{align}
Moreover, we assume that the speed of the delay satisfies
\begin{align}
\tau'(t)\leq d<1, \quad \forall\ t>0, \label{1.3}
\end{align}
that
\begin{align}
\tau\in W^{2,\infty}([0,T]), \quad \forall\ T>0 \label{1.4}
\end{align}
and that $a_{0}$, $a_{1}$ satisfy
\begin{align}
|a_{1}|<\sqrt{1-d}\, a_{0}. \label{1.5}
\end{align}

As in \cite{np2011}, let us introduce the function
\begin{equation}\label{2.5}
z(x,\rho ,t)=u_t(x, t-\tau(t)\rho ),\quad x\in\Omega,\ \rho\in (0,1), \ t>0.
\end{equation}
Then, problem \eqref{1.1} is equivalent to
\bes\left\{\begin{array}{ll}
\displaystyle  u_{tt}(x,t)-\Delta u(x,t)+\int_{0}^{t}g(t-s)\Delta
u(x,s)\, {\rm d}s\quad &  \\\medskip
\quad \quad \quad \quad \quad \quad +a_{0}u_{t}(x,t)+a_{1}z(x,1 ,t)=0,
&(x,t)\in \Omega \times (0,\infty ), \\\medskip
\displaystyle \tau(t) z_t(x,\rho,t)+(1-\tau'(t)\rho)z_{\rho}(x,\rho, t)=0, &(x,\rho,t)\in \Omega\times (0,1) \times (0,\infty ),\\\medskip
  \displaystyle u(x,t)=0,  &(x,t)\in \partial\Omega \times [0,\infty ),\\\medskip
  \displaystyle  z(x,0,t)=u_t(x,t), &(x,t)\in \Omega \times (0,\infty ), \\\medskip
\displaystyle u(x,0)=u_{0}(x),\quad u_{t}(x,0)=u_{1}(x), \quad   &x
\in  \Omega,\\
  z(x,\rho ,0)=f_0(x,-\rho \tau(0) ),\quad  &(x,\rho)\in \Omega\times (0,1),
\end{array}\right.\label{2.6}
 \ees

 We now state, without a proof, a well-posedness result,
which can be established by combining the arguments
of \cite{fnv2010,ks2011}.

\begin{lemma}  \label{th2.1}
Let \eqref{1.2}--\eqref{1.5} be satisfied and $g$ satisfies {\rm (G1)}. Then given $u_0\in H^1_0 (\Omega)$, $u_1\in L^2 (\Omega)$,
$f_0\in L^2(\Omega\times (0,1))$ and $T>0$, there exists a unique weak solution $(u, z)$ of the problem \eqref{2.6} on $(0, T)$ such that
$$u\in C(0,T; H^1_0 (\Omega))\cap C^1(0,T; L^2 (\Omega)),\quad u_t\in L^2(0,T; H^1_0 (\Omega))\cap L^2((0,T)\times\Omega).$$
\end{lemma}

Inspired by \cite{m2008,np2011}, we define the new energy functional as
\begin{align}
E(t):=\frac{1}{2}\int_{\Omega}\left[u_{t}^{2}+\left(1-\int_{0}^{t}g(s){\rm
d}s\right)|\nabla u|^{2}\right]{\rm d}x+\frac{1}{2}(g\circ\nabla
u)(t)\nonumber\\+\frac{\xi}{2}\int_{t-\tau(t)}^{t}\int_{\Omega}e^{\lambda(s-t)}u_{t}^{2}(x,s){\rm
d}x{\rm d}s, \label{3.1}
\end{align}
where $\xi$, $\lambda$ are suitable positive constants, and
\begin{align*}
(g\circ v)(t)=\int_{\Omega}\int_{0}^{t}g(t-s)|v(t)-v(s)|^{2}{\rm
d}s{\rm d}x, \quad \forall\ v\in L^{2}(\Omega).
\end{align*}
We will fix $\xi$ such that
\begin{align}
2a_{0}-\frac{|a_{1}|}{\sqrt{1-d}}-\xi>0 \quad \text{and} \quad
\xi-\frac{|a_{1}|}{\sqrt{1-d}}>0, \label{3.2}
\end{align}
and
\begin{align}
\lambda<\frac{1}{\overline{\tau}}\left|\log\frac{|a_{1}|}{\xi\sqrt{1-d}}\right|.
\label{3.3}
\end{align}
In fact, the existence of such a constant $\xi$ is guaranteed by the
assumption \eqref{1.5}.

Our main result is the following.

\begin{theorem}  \label{th2.2}
Let \eqref{1.2}--\eqref{1.5} be satisfied and $g$ satisfies {\rm (G1)} and {\rm (G2)}. Then there exist two positive constants $K, k$ such that, for any
solution of problem \eqref{1.1}, the energy  satisfies
\begin{align}
E(t)\leq Ke^{-k\int_{t_{0}}^{t}\xi(s){\rm d}s}, \quad \forall\ t\geq
t_{0}.
\label{2.10}
\end{align}
\end{theorem}

\begin{remark} \label{re1}
Note that the exponential or polynomial decay estimate  is
only a particular case of \eqref{2.10}. More precisely, we can
obtain exponential decay for $\xi(t)\equiv a$ and polynomial decay
for $\xi(t)\equiv a(1+t)^{-1}$, where $a > 0$ is a constant.
\end{remark}

\begin{remark} \label{re2}
Estimate \eqref{2.10} is also true for $t\in [0, t_0]$ by virtue of
the continuity and boundedness of $E(t)$ and $\xi(t)$.
\end{remark}

\section{General decay of the solution} \label{s3}

As mentioned earlier, the proof of the general decay result is given in this section. We have the following lemmas.

\begin{lemma}  \label{pr3.1}
Let \eqref{1.2}--\eqref{1.5} be satisfied  and $g$ satisfies {\rm (G1)}. Then for all regular
solution of problem \eqref{1.1}, the energy functional defined by \eqref{3.1} is  non-increasing and
satisfies
\begin{align}
E'(t)\leq&\frac{1}{2}(g'\circ\nabla
u)(t)-\frac{1}{2}g(t)\int_{\Omega}|\nabla u|^{2}{\rm
d}x-C_{1}\int_{\Omega}\left[u_{t}^{2}(x,t)+u_{t}^{2}(x,t-\tau(t))\right]{\rm
d}x \nonumber \\
&-\frac{\lambda\xi}{2}\int_{t-\tau(t)}^{t}\int_{\Omega}e^{\lambda(s-t)}u_{t}^{2}(x,s){\rm
d}x{\rm d}s\le 0 \label{3.4}
\end{align}
for some positive constant $C_{1}$.
\end{lemma}

{\bf Proof.}\ Differentiating \eqref{3.1} and by \eqref{1.1}, we obtain
\begin{align*}
E'(t)=&\int_{\Omega}\left[u_{t}u_{tt}+\left(1-\int_{0}^{t}g(s){\rm d}s\right)\nabla
u\cdot\nabla u_{t}-\frac{1}{2}g(t)|\nabla u|^{2}\right]{\rm d}x
\\
&+\int_{0}^{t}g(t-s)\int_{\Omega}\nabla u_{t}(t)\cdot[\nabla
u(t)-\nabla u(s)]{\rm d}x{\rm d}s +\frac{1}{2}\int_{0}^{t}g'(t-s)\int_{\Omega}|\nabla u(t)-\nabla
u(s)|^{2}{\rm d}x{\rm d}s \\
&+\frac{\xi}{2}\int_{\Omega}u_{t}^{2}(x,t){\rm
d}x-\frac{\xi}{2}\int_{\Omega}e^{-\lambda\tau(t)}u_{t}^{2}(x,t-\tau(t))(1-\tau'(t)){\rm
d}x \\
&-\lambda\frac{\xi}{2}\int_{t-\tau(t)}^{t}\int_{\Omega}e^{-\lambda(t-s)}u_{t}^{2}(x,s){\rm
d}x{\rm d}s \\
=&\int_{\Omega}\left[u_{t}u_{tt}+\nabla u\cdot\nabla
u_{t}-\int_{0}^{t}g(t-s)\nabla u(s)\cdot\nabla u_{t}(t){\rm d}s
\right]{\rm d}x \\
&-\frac{1}{2}g(t)\int_{\Omega}|\nabla u|^{2}{\rm
d}x+\frac{1}{2}(g'\circ\nabla
u)(t)+\frac{\xi}{2}\int_{\Omega}u_{t}^{2}(x,t){\rm d}x \\
&-\frac{\xi}{2}\int_{\Omega}e^{-\lambda\tau(t)}u_{t}^{2}(x,t-\tau(t))(1-\tau'(t)){\rm
d}x  -\lambda\frac{\xi}{2}\int_{t-\tau(t)}^{t}\int_{\Omega}e^{-\lambda(t-s)}u_{t}^{2}(x,s){\rm
d}x{\rm d}s,
\end{align*}
and then, using integration by parts, the assumptions
\eqref{1.2}--\eqref{1.3} and some manipulations as in \cite{np2011},
\begin{align}
E'(t)=&-a_{0}\int_{\Omega}u_{t}^{2}(x,t){\rm
d}x-a_{1}\int_{\Omega}u_{t}(t)\int_{\Omega}u_{t}(t-\tau(t)){\rm d}x-\frac{1}{2}g(t)\int_{\Omega}|\nabla u|^{2}{\rm
d}x
\nonumber \\
&+\frac{1}{2}(g'\circ\nabla
u)(t)+\frac{\xi}{2}\int_{\Omega}u_{t}^{2}(x,t){\rm d}x -\frac{\xi}{2}\int_{\Omega}e^{-\lambda\tau(t)}u_{t}^{2}(x,t-\tau(t))(1-\tau'(t)){\rm
d}x \nonumber \\
&-\lambda\frac{\xi}{2}\int_{t-\tau(t)}^{t}\int_{\Omega}e^{-\lambda(t-s)}u_{t}^{2}(x,s){\rm
d}x{\rm d}s \nonumber \\
\leq&-a_{0}\int_{\Omega}u_{t}^{2}(x,t){\rm
d}x-a_{1}\int_{\Omega}u_{t}(t)\int_{\Omega}u_{t}(t-\tau(t)){\rm d}x
-\frac{1}{2}g(t)\int_{\Omega}|\nabla u|^{2}{\rm
d}x\nonumber \\
&+\frac{1}{2}(g'\circ\nabla
u)(t)+\frac{\xi}{2}\int_{\Omega}u_{t}^{2}(x,t){\rm d}x-\frac{\xi}{2}(1-d)e^{-\lambda\overline{\tau}}\int_{\Omega}u_{t}^{2}(x,t-\tau(t)){\rm
d}x \nonumber \\
&-\lambda\frac{\xi}{2}\int_{t-\tau(t)}^{t}\int_{\Omega}e^{-\lambda(t-s)}u_{t}^{2}(x,s){\rm
d}x{\rm d}s \nonumber \\
\leq&\frac{1}{2}(g'\circ\nabla
u)(t)-\frac{1}{2}g(t)\int_{\Omega}|\nabla u|^{2}{\rm
d}x-\left(a_{0}-\frac{|a_{1}|}{2\sqrt{1-d}}-\frac{\xi}{2}\right)\int_{\Omega}u_{t}^{2}(x,t){\rm
d}x \nonumber \\
&-\left(e^{-\lambda\overline{\tau}}\frac{\xi}{2}(1-d)-\frac{|a_{1}|}{2}\sqrt{1-d}\right)\int_{\Omega}u_{t}^{2}(x,t-\tau(t)){\rm
d}x \nonumber \\
&-\lambda\frac{\xi}{2}\int_{t-\tau(t)}^{t}\int_{\Omega}e^{-\lambda(t-s)}u_{t}^{2}(x,s){\rm
d}x{\rm d}s. \label{3.5}
\end{align}
Combining \eqref{3.2}--\eqref{3.3}, \eqref{3.5} and hypothese
{\rm (G1)}, \eqref{3.4} is established. \quad $\Box$

Now we are going to construct a Lyapunov functional $L$ equivalent
to $E$. For this purpose, we define the following functionals:
\begin{align}
I(t):=\int_{\Omega}uu_{t}{\rm d}x, \label{3.6}
\end{align}
\begin{align}
K(t):=-\int_{\Omega}u_{t}\int_{0}^{t}g(t-s)(u(t)-u(s)){\rm d}s{\rm
d}x, \label{3.7}
\end{align}

Set
\begin{align}
L(t)=NE(t)+\varepsilon I(t)+K(t) \label{3.8}
\end{align}
where $N$ and $\varepsilon$ are suitable positive constants to be
determined later. Similar as in \cite{m2008}, we can prove that, for
$\varepsilon$ small enough while $N$ large enough, there exist two
positive constants $\beta_{1}$, $\beta_{2}$ such that
\begin{align}
\beta_{1}E(t)\leq L(t)\leq\beta_{2}E(t), \quad \forall\ t\geq0.
\label{3.9}
\end{align}

The following estimates hold true.
\begin{lemma} \label{le3.2}
Under the assumption {\rm (G1)}, the functional $I$ satisfies, along the
solution, the estimate
\begin{align}
I'(t)\leq-\frac{l}{2}\int_{\Omega}|\nabla u|^{2}{\rm
d}x+C_{2}\int_{\Omega}[u_{t}^{2}(x,t)+u_{t}^{2}(x,t-\tau(t))]{\rm
d}x+C_{3}(g\circ\nabla u)(t). \label{3.10}
\end{align}
\end{lemma}

{\bf Proof.}\ Differentiating and integrating by parts
\begin{align}
I'(t)=&\int_{\Omega}u_{t}^{2}{\rm d}x+\int_{\Omega}u\left(\Delta
u-\int_{0}^{t}g(t-s)\Delta u(s){\rm
d}s-a_{0}u_{t}(t)-a_{1}u_{t}(t-\tau(t))\right){\rm d}x \nonumber \\
=&\int_{\Omega}u_{t}^{2}{\rm d}x-l\int_{\Omega}|\nabla u|^{2}{\rm
d}x+\int_{\Omega}\nabla u\cdot\int_{0}^{t}g(t-s)(\nabla u(s)-\nabla
u(t)){\rm d}s{\rm d}x \nonumber
\\&-a_{0}\int_{\Omega}u(t)u_{t}(t){\rm
d}x-a_{1}\int_{\Omega}u(t)u_{t}(t-\tau(t)){\rm d}x. \label{3.11}
\end{align}
Now, using Young's inequality and {\rm (G1)}, we obtain (see \cite{m2008})
\begin{align}
&\int_{\Omega}\nabla u\cdot\int_{0}^{t}g(t-s)(\nabla u(s)-\nabla
u(t)){\rm d}s{\rm d}x \nonumber \\
\leq&\delta\int_{\Omega}|\nabla u|^{2}{\rm d}x+\frac{1}{4\delta}
\int_{\Omega}\left(\int_{0}^{t}g(t-s)|\nabla u(s)-\nabla
u(t)|{\rm d}s\right)^{2}{\rm d}x \nonumber \\
\leq&\delta\int_{\Omega}|\nabla u|^{2}{\rm
d}x+\frac{1-l}{4\delta}(g\circ\nabla u)(t),\quad \forall\
\delta>0.\label{3.12}
\end{align}
Also, using Young's and Poincar\'e's inequalities gives
\begin{align}
-a_{0}\int_{\Omega}u(t)u_{t}(t){\rm
d}x\leq\delta\int_{\Omega}|\nabla u|^{2}{\rm
d}x+C(\delta)\int_{\Omega}u_{t}^{2}{\rm d}x, \label{3.13}
\end{align}
\begin{align}
-a_{1}\int_{\Omega}u(t)u_{t}(t-\tau(t)){\rm
d}x\leq\delta\int_{\Omega}|\nabla u|^{2}{\rm
d}x+C(\delta)\int_{\Omega}u_{t}^{2}(t-\tau(t)){\rm d}x. \label{3.14}
\end{align}
Combining \eqref{3.11}--\eqref{3.14} and choosing $\delta$ small
enough, we obtain \eqref{3.10}. \quad $\Box$

\begin{lemma} \label{le3.3}
Under the assumption {\rm (G1)}, the functional $K$ satisfies, along the
solution, the estimate
\begin{align}
K'(t)\leq&-\left(\int_{0}^{t}g(s){\rm
d}s-2\delta\right)\int_{\Omega}u_{t}^{2}{\rm
d}x+\delta\int_{\Omega}|\nabla u|^{2}{\rm
d}x+\frac{C_{4}}{\delta}(g\circ\nabla u)(t) \nonumber \\
&-\frac{C_{5}}{\delta}(g'\circ\nabla
u)(t)+\delta\int_{\Omega}u_{t}^{2}(t-\tau(t)){\rm d}x. \label{3.15}
\end{align}
\end{lemma}

{\bf Proof.}\ By exploiting \eqref{1.1} and integrating by parts, we
have
\begin{align*}
K'(t)=&\left(1-\int_{0}^{t}g(s){\rm d}s\right)\int_{\Omega}\nabla
u\cdot\int_{0}^{t}g(t-s)(\nabla u(t)-\nabla
u(s)){\rm d}s{\rm d}x  \\
&+\int_{\Omega}\left(\int_{0}^{t}g(t-s)(\nabla u(s)-\nabla
u(t)){\rm d}s\right)^{2}{\rm d}x-\int_{\Omega}u_{t}\int_{0}^{t}g'(t-s)(u(t)-u(s)){\rm d}s{\rm
d}x  \\
&-\int_{0}^{t}g(s){\rm d}s\int_{\Omega}u_{t}^{2}{\rm d}x   +\int_{\Omega}\left(\int_{0}^{t}g(t-s)(u(t)-u(s)){\rm
d}s\right)[a_{0}u_{t}(t)+a_{1}u_{t}(t-\tau(t))]{\rm d}x.
\end{align*}
Using Young's and Poincar\'e's inequalities, we obtain (see
\cite{m2008,ly2011})
\begin{align*}
&\left(1-\int_{0}^{t}g(s){\rm d}s\right)\int_{\Omega}\nabla
u\cdot\int_{0}^{t}g(t-s)(\nabla u(t)-\nabla u(s)){\rm d}s{\rm d}x \\
\leq&\delta\int_{\Omega}|\nabla u|^{2}{\rm
d}x+\frac{C}{\delta}(g\circ\nabla u)(t),
\end{align*}
\begin{align*}
-\int_{\Omega}u_{t}\int_{0}^{t}g'(t-s)(u(t)-u(s)){\rm d}s{\rm
d}x\leq\delta\int_{\Omega}u_{t}^{2}{\rm
d}x-\frac{C}{\delta}(g'\circ\nabla u)(t),
\end{align*}
\begin{align*}
&\int_{\Omega}\left(\int_{0}^{t}g(t-s)(u(t)-u(s)){\rm
d}s\right)[a_{0}u_{t}(t)+a_{1}u_{t}(t-\tau(t))]{\rm d}x \\
\leq&\frac{C}{\delta}(g\circ\nabla
u)(t)+\delta\int_{\Omega}u_{t}^{2}{\rm
d}x+\delta\int_{\Omega}u_{t}^{2}(t-\tau(t)){\rm d}x.
\end{align*}
Combining all above estimates, \eqref{3.15} is established. \quad
$\Box$

Now, we are ready to prove the general decay result.

{\bf Proof of Theorem \ref{th2.2}.} \ Since $g$ is positive, we
have, for any $t_{0}>0$,
\begin{align*}
\int_{0}^{t}g(s){\rm d}s\geq\int_{0}^{t_{0}}g(s){\rm d}s:=g_{0}>0,
\quad t\geq t_{0}.
\end{align*}
By using \eqref{3.4}, \eqref{3.8}, \eqref{3.10} and \eqref{3.15}, a
series of computations yields, for $t\geq t_{0}$,
\begin{align}
L'(t)\leq&\frac{N}{2}(g'\circ\nabla
u)(t)-\frac{N}{2}g(t)\int_{\Omega}|\nabla u|^{2}{\rm
d}x-NC_{1}\int_{\Omega}[u_{t}^{2}(x,t)+u_{t}^{2}(x,t-\tau(t))]{\rm
d}x \nonumber \\
&-\frac{\lambda\xi
N}{2}\int_{t-\tau(t)}^{t}\int_{\Omega}e^{\lambda(s-t)}u_{t}^{2}(x,s){\rm
d}x{\rm d}s+\varepsilon
C_{2}\int_{\Omega}[u_{t}^{2}(x,t)+u_{t}^{2}(x,t-\tau(t))]{\rm
d}x \nonumber \\
&-\frac{\varepsilon l}{2}\int_{\Omega}|\nabla u|^{2}{\rm
d}x+\varepsilon C_{3}(g\circ\nabla u)(t)-\left(\int_{0}^{t}g(s){\rm
d}s-2\delta\right)\int_{\Omega}u_{t}^{2}{\rm
d}x+\delta\int_{\Omega}|\nabla u|^{2}{\rm d}x \nonumber \\
&+\frac{C_{4}}{\delta}(g\circ\nabla u)(t)
-\frac{C_{5}}{\delta}(g'\circ\nabla
u)(t)+\delta\int_{\Omega}u_{t}^{2}(t-\tau(t)){\rm d}x \nonumber \\
=&-[(NC_{1}+g_{0})-2\delta-\varepsilon
C_{2}]\int_{\Omega}u_{t}^{2}{\rm d}x+\left(\varepsilon
C_{3}+\frac{1}{\delta}C_{4}\right)(g\circ\nabla
u)(t)\nonumber\\
&+\left(\frac{N}{2}-\frac{C_{5}}{\delta}\right)(g'\circ\nabla
u)(t)-\left(\frac{\varepsilon l}{2}-\delta\right)\int_{\Omega}|\nabla
u|^{2}{\rm d}x \nonumber \\
&-(NC_{1}-\delta-\varepsilon
C_{2})\int_{\Omega}u_{t}^{2}(x,t-\tau(t)){\rm
d}x \nonumber \\
&-\frac{\lambda\xi
N}{2}\int_{t-\tau(t)}^{t}\int_{\Omega}e^{\lambda(s-t)}u_{t}^{2}(x,s){\rm
d}x{\rm d}s. \label{3.16}
\end{align}
At this point, we choose $\varepsilon$ small enough such that
$\varepsilon<\frac{g_{0}}{2C_{2}}$ and \eqref{3.9} hold, and
$\delta$ sufficiently small such that
$$\alpha_{1}=\frac{\varepsilon l}{2}-\delta>0.$$
As long as $\varepsilon$ and $\delta$ are fixed, we choose $N$ large
enough such that
\begin{align*}
NC_{1}-2\delta>0, \quad \alpha_{2}=NC_{1}-\delta-\varepsilon C_{2}>0
\quad \text{and} \quad
\alpha_{3}=\frac{N}{2}-\frac{C_{5}}{\delta}>0.
\end{align*}
Thus, it follows from {\rm (G2)} and \eqref{3.16} that
\begin{align}
L'(t)\leq&\frac{g_{0}}{2}\int_{\Omega}u_{t}^{2}{\rm
d}x-\alpha_{1}\int_{\Omega}|\nabla u|^{2}{\rm d}x
-\alpha_{2}\int_{\Omega}u_{t}^{2}(x,t-\tau(t)){\rm d}x\nonumber
\\&-\frac{\lambda\xi
N}{2}\int_{t-\tau(t)}^{t}\int_{\Omega}e^{\lambda(s-t)}u_{t}^{2}(x,s){\rm
d}x{\rm d}s+\alpha_{4}(g\circ\nabla u)(t) \nonumber\\
\leq&-C_{6}E(t)+\alpha_{4}(g\circ\nabla u)(t), \quad \forall\ t\geq
t_{0}, \label{3.17}
\end{align}
where $\alpha_{4}=\varepsilon C_{3}+\frac{1}{\delta}C_{4}>0$. It
follows from \eqref{3.17}, {\rm (G1)} and \eqref{3.4} that
\begin{align}
\xi(t)L'(t)\leq&-C_{6}\xi(t)E(t)+\alpha_{4}\xi(t)(g\circ\nabla u)(t) \nonumber\\
\leq&-C_{6}\xi(t)E(t)-\alpha_{4}(g'\circ\nabla u)(t) \nonumber\\
\leq&-C_{6}\xi(t)E(t)-C_{7}E'(t),\quad \forall\ t\geq t_{0},
\label{3.18}
\end{align}
That is
\begin{align}
F'(t)\leq-C_{8}\xi(t)E(t)\leq-k\xi(t)F(t), \quad \forall\ t\geq
t_{0}, \label{3.19}
\end{align}
where $$F(t)=\xi(t)L(t)+C_{7}E(t)$$ is clearly equivalent to $E(t)$
and $k$ is a positive constant.

Consequently, \eqref{2.10} can be obtained by \eqref{3.9} and \eqref{3.19}.\quad $\Box$

\subsection*{Acknowledgments}
This work was partly supported by the Tianyuan Fund of Mathematics (Grant No. 11026211) and the Natural
Science Foundation of the Jiangsu Higher Education Institutions (Grant No. 09KJB110005).

\end{document}